\documentclass[twoside,12pt]{article}
\usepackage{amsmath,amssymb,epsfig}
\date{}
\newtheorem{Lemma}{LEMMA}[section]
\newtheorem{Corollary}[Lemma]{COROLLARY}
\newtheorem{Theorem}[Lemma]{THEOREM}
\newtheorem{Proposition}[Lemma]{PROPOSITION}

\newcommand{\bnum}{\begin{enumerate}}
\newcommand{\enum}{\end{enumerate}}
\newcommand{\bi}{\begin{itemize}}
\newcommand{\ei}{\end{itemize}}
\newcommand{\btab}{\begin{tabular}}
\newcommand{\etab}{\end{tabular}}
\newcommand{\beq}{\begin{eqnarray*}}
\newcommand{\eeq}{\end{eqnarray*}}
\newcommand{\beqn}{\begin{eqnarray}}
\newcommand{\eeqn}{\end{eqnarray}}

\newcommand{\bq}{\begin{equation}}
\newcommand{\eq}{\end{equation}}

\def\phi{\varphi}
\def\epsilon{\varepsilon}

\newcommand{\BR}{\mathbb R}
\newcommand{\BC}{\mathbb C}
\newcommand{\BH}{\mathbb H}

\newcommand{\kasten}{\vbox{\hrule height 8pt width 8.6pt depth -7.4pt
    \hbox{\vrule width 0.6pt height 7.4pt
    \kern 7.4pt \vrule width 0.6pt height 7.4pt}
    \hrule height 0.6pt width 8.6pt}}
\newcommand{\ok}{\hfill\kasten}
\newcommand{\bpf}{\begin{Proof}}
\newcommand{\epf}{\ok\end{Proof}\bigskip\noindent}
\newcommand{\bthm}{\begin{Theorem}}
\newcommand{\ethm}{\end{Theorem}}
\newcommand{\ble}{\begin{Lemma}}
\newcommand{\ele}{\end{Lemma}}
\newcommand{\bprop}{\begin{Proposition}}
\newcommand{\eprop}{\end{Proposition}}
\newcommand{\bcor}{\begin{Corollary}}
\newcommand{\ecor}{\end{Corollary}}
\begin{document}
\title{A characterization of Clifford parallelism by automorphisms}

\author{Rainer L\"owen}

\maketitle

\thispagestyle{empty}

\begin{abstract}
Betten and Riesinger have shown that Clifford parallelism on real projective space is the only topological parallelism that is left 
invariant by a group of dimension at least 5. We improve the bound to 4. Examples of different parallelisms admitting a group of 
dimension $\le 3$ are known, so 3 is the `critical dimension'.
 
MSC 2000: 51H10, 51A15, 51M30 
\end{abstract}

Consider $\BR^4$ as the quaternion skew field $\BH$. Then the orthogonal group $\mathop{\rm SO}(4,\BR)$ may be described as the product of two 
commuting copies $\tilde \Lambda, \tilde \Phi$ of the unitary group $\mathop{\rm U}(2,\BC)$, consisting 
of the maps $q \mapsto aq$ and $q\mapsto qb$, respectively, 
where $a$, $b$ are 
quaternions of norm one and multiplication is quaternion multiplication. The intersection of the two factors is of order two, 
containing the map $-\rm id$.
Thus, passing to projective space, we get $\mathop{\rm PSO}(4,\BR) = \Lambda \times \Phi$, 
a direct product of two copies of $\rm SO (3,\BR)$.
The left and right Clifford parallelisms are defined as the equivalence relations on the line space of $\rm PG(3,\BR)$ formed by the 
orbits of $\Lambda$ and $\Phi$, respectively. 

The two Clifford parallelisms are equivalent under quaternion conjugation $q \to \bar q$; 
this is immediate from their definition in view of the fact that conjugation does not change the norm and is an anti-automorphism, i.e., 
that $\overline {pq} = \bar q \bar p$. Note that both $\Lambda$ and $\Phi$ are transitive on the point set of projective space. Since they centralize
one another, each acts transitively on the parallelism defined by the other, and the group $\mathop{\rm PSO}(4,\BR)$ leaves both parallelisms invariant
(we say that it consists of \it automorphisms \rm of these parallelisms). For more information on Clifford parallels, see 
\cite{berg}, \cite{kling}, \cite{cliff}. For generalizations to other dimensions, compare also \cite{tyr}.

The notion of a \it topological parallelism \rm on real projective 3-space $\mathop {\rm PG}(3,\BR)$ generalizes this example. 
A \it spread \rm is a set $\cal C$ of lines such that every point is incident with exactly one of them, and
a topological parallelism may be defined as a compact set $\Pi$ of compact spreads such that every line belongs to exactly one of them; 
see, e.g., \cite{coll} for details.
Many examples of different topological parallelisms have been constructed in a series of papers by Betten and Riesinger, see, e.g., \cite{gl} 
and references given therein. 

The group $\Sigma = \mathop{\rm Aut} \Pi$ of automorphisms of a topological parallelism is a closed subgroup of 
$\mathop {\rm PGL}(4,\BR)$. The identity component $\Delta = \Sigma^1$ is compact \cite{comp}, and hence (conjugate to) a subgroup of
$\mathop{\rm PSO}(4,\BR) \cong \mathop{\rm SO}(3,\BR) \times \mathop{\rm SO}(3,\BR)$.
Hence, in the case of the Clifford parallelism, $\Delta$ is the 6-dimensional group 
$\mathop{\rm PSO}(4,\BR)$ that we used to define the parallelism. 
Betten and Riesinger \cite{coll} proved that no 
other topological parallelism has a group of dimension $\dim \Sigma \ge 5$. 
Examples of parallelisms with 1-, 2- or 3-dimensional automorphism groups are known, see \cite{reg}, \cite{gl}, \cite{nonreg}.
Here we consider parallelisms with a 4-dimensional group.

\bthm
Let $\Sigma$ be the automorphism group of a topological parallelism $\Pi$ on $\rm PG(3,\BR)$. If 
$\dim \Sigma \ge 4$, then $\Pi$ is equivalent to the Clifford parallelism. 
\ethm

\bpf
Recall that a topological parallelism $\Pi$ is 
homeomorphic to the real projective plane in the Hausdorff topology on the space of compact sets of lines, and that every equivalence class is a 
compact spread and homeomorphic to the 2-sphere, compare \cite{coll}.

It suffices to consider the identity component $\Delta = \Sigma^1$, and
we may assume that $ \dim \Delta = 4$. Further, up to equivalence, we may assume that $\Delta = \Lambda \cdot \Gamma$, where 
$\Gamma \le \Phi$ is the subgroup defined by restricting the factor $b$ to be a complex number (here we use the notation of the introduction.) 
Since $\Lambda$ does not have any one-dimensional coset spaces, we know that $\Lambda$ 
acts on $\Pi$ either transitively or trivially. If it acts trivially, then the classes of $\Pi$ are the $\Lambda$-orbits of lines, and we 
have the Clifford parallelism. Observe here that every $\Lambda$-orbit is contained in a single class, and both the orbit and the class are 2-spheres.

In what follows, assume therefore that $\Lambda$ acts transitively on $\Pi$. There is only one possibility for this action namely, the 
standard transitive action of $\mathop{\rm SO}(3,\BR)$ on the real projective plane.
Every 2-dimensional subgroup of 
$\Delta$ contains $\Gamma$. Hence, there is no effective action of $\Delta$ 
on the projective plane $\Pi$, and the kernel can only be $\Gamma$ since the only other proper normal subgroup is $\Lambda$, which is transitive. 
If ${\cal C} \in \Pi$ is any equivalence class, then the stabilizer $\Lambda_{\cal C}$ is a product of a 1-torus and a group of order two. 
Hence $\Delta_{\cal C}$ contains a 2-torus $T$. There is only one conjugacy class of 2-tori in $\Delta$, represented by the group
 $$T_0 = \{\langle q \rangle \mapsto \langle aqb \rangle \ \vert \ a,b\in \BC, \vert a\vert = \vert b \vert =1\}.$$ 
Here, $\langle q \rangle$ denotes the one-dimensional real vector space spanned by $q$. 
We may assume that $T = T_0$. Write quaternions as pairs of complex numbers with multiplication $(x,y)(u,v) = (xu-\bar v y, vx + y\bar u)$, see
\cite{CPP}, 11.1. Then complex numbers become pairs $(a,0)$, and the elements of $T$ are now given by
  $$\langle (z,w) \rangle \mapsto \langle (azb,aw\bar b)\rangle .$$
The kernel of ineffectivity 
of $T$ on the 2-sphere $\cal C$ must be a 1-torus $\Xi$, and the elements of the kernel
other than the identity cannot have eigenvalue 1 --- otherwise they would be axial collineations 
of the translation plane defined by the spread $\cal C$ 
and would act non-trivially on $\cal C$. There are only two subgroups of the 2-torus satisfying these conditions, 
given by $b = 1$ and by $a = 1$, respectively. In other words, the kernel $\Xi$ is a subgroup either of $\Lambda$ or of $\Phi$. In both 
cases, $\cal C$ consists of the fixed lines of $\Xi$. If $\Xi\le \Phi$, then $\Lambda$ permutes these lines, contrary to the 
transitivity of $\Lambda$ on $\Pi$. If $\Xi \le \Lambda$, then $\Phi$ permutes the fixed lines, which means that $\cal C$ is a $\Phi$-orbit. 
Now $\Lambda$ is transitive both on $\Pi$ and on the set of $\Phi$-orbits, hence $\Pi$ equals the Clifford parallelism formed by the $\Phi$-orbits.

\epf


\bibliographystyle{plain}

\bigskip
\bigskip
\noindent{Rainer L\"owen, Institut f\"ur Analysis und Algebra,
Technische Universit\"at Braunschweig,
Pockelsstra{\ss}e 14,
D 38106 Braunschweig,
Germany}

\end{document}